\documentclass[oneside,11pt]{amsart}
\usepackage{amsmath}
\usepackage{amssymb}
\usepackage{graphicx}
\usepackage{tikz}
\usepackage[colorlinks=false]{hyperref}
\usepackage{enumitem}
 \usepackage{float} 
 \usepackage{caption}
\usepackage{color}
\usepackage{soul}

\newtheorem{thm}{Theorem}[section]
\newtheorem{prop}[thm]{Proposition}
\newtheorem{lem}[thm]{Lemma}

\theoremstyle{definition}
\newtheorem{defn}[thm]{Definition}
\newtheorem{rem}[thm]{Remark}

\renewcommand{\emptyset}{\varnothing}

\renewcommand{\setminus}{-}

\newcommand{\field}[1]{\mathbb{#1}}
\newcommand{\Z}{\field{Z}}

\renewcommand{\implies}{\Rightarrow}

\DeclareMathOperator{\diam}{diam}

\usepackage{ifthen}

\newcommand{\showcomments}{yes}

\newsavebox{\commentbox}
{\ifthenelse{\equal{\showcomments}{yes}}%
{\footnotemark
    \begin{lrbox}{\commentbox}
    \begin{minipage}[t]{1.25in}\raggedright\sffamily\upshape\tiny
    \footnotemark[\arabic{footnote}]}%
{\begin{lrbox}{\commentbox}}}%
{\ifthenelse{\equal{\showcomments}{yes}}%
{\end{minipage}\end{lrbox}\marginpar{\usebox{\commentbox}}}%
{\end{lrbox}}}

\newcounter{acomments}

\newcounter{ccomments}

\newcounter{hcomments}

\author{Hoang Thanh Nguyen}
\address{Department of Mathematics\\
The University of Danang - University of Science and
Education\\
 459 Ton Duc Thang, Da Nang, Vietnam}
\email{nthoang.math@gmail.com}

\author{Hung Cong Tran}
\address{University of Oklahoma, Norman, OK 73019-3103, USA}
\email{Hung.C.Tran-1@ou.edu}

\title{ Strongly quasiconvex subgroups in graphs of groups}

\begin{document}

\date{\today}
\maketitle

\begin{abstract}
Given a graph of groups $\mathcal{G} = (\Gamma, \{G_v\}, \{G_e\})$ with certain conditions on vertex groups and $G$ acts acylindrically on its Bass-Serre tree $T$. Let $H$ be a  finitely generated subgroup  of $G$. We prove the following statements equivalence.
\begin{enumerate}
    \item $H$ has finite height
    \item $(G, T, H)$ is a $A/QI$--triple
    \item $H$ is strongly quasiconvex and virtually free in $G$.
\end{enumerate}

We also give a condition to determine whether strong quasiconvexity in a group is preserved under amalgams. 
\end{abstract}

\section{Introduction
}
The {\it height} of a finitely generated subgroup $H$ in a finitely generated group $G$ is the minimal $n \in \mathbb N_{\ge 0}$ such that for any distinct cosets $g_0 H,\dots,g_n H \in G/H$, the intersection
$g_0 H g_0^{-1} \cap \dots \cap g_n H g_n^{-1}$
is finite (the height is \textit{infinite} if no such $n$ exists).
The subgroup $H$ is called {\it strongly quasi-convex} in $G$ if for any $L \ge 1$, $C \ge 0$ there exists $R = R(L,C)$ such that every $(L,C)$--quasi-geodesic in $G$ with endpoints in $H$ is contained in the $R$--neighborhood of $H$ (see \cite{Tra19} or \cite{Gen19}). We note that strong quasiconvexity does not depend on the choice of a finite generating set of the ambient group and it agrees with quasiconvexity when the ambient group is hyperbolic.
Introduced recently by Abbott-Manning \cite{AM21}, an $A/QI$--triple $(G, X, H)$ consists of a Gromov hyperbolic space
$X$ and an action $G \curvearrowright X$ which is {\it acylindrical along} the finitely generated infinite
subgroup $H$, and so that $h \to hx$ (for some $x \in X$) gives a quasi-isometric embedding of $H$ into $X$.

In \cite{GMRS98}, the authors prove that strongly quasi-convex subgroups in hyperbolic groups have finite height. It is a long-standing question asked by Swarup that whether or not the converse is
true (see Question 1.8 in \cite{Bes}). If the converse is true, then we could characterize strongly quasi-convex subgroup of a finitely generated group purely in terms
of group theoretic notions.
Tran in \cite{Tra19} generalizes the result of \cite{GMRS98} by showing that strongly quasi-convex subgroups in any finitely generated group have finite height. Abott-Manning in \cite{AM21} show that if $(G, X,H)$ is an $A/QI$--triple then $H$ has finite height in $G$. It is natural and reasonable to extend the Swarup’s question to
strongly quasiconvex subgroups of finitely generated groups to ask whether or not the converse is true, i.e. finite height implies strong quasiconvexity, $A/QI$--triple. 

A counter example to the above question has been addressed in \cite{NTY21} where the authors show that the fundamental group of a torus bundle $M$ over circle with Anosov monodromy contains finite height subgroups that are not strongly quasiconvex. We note that the action of fundamental group $\pi_1(M)$ on the tree $T$ is not acylindrical. In the setting of 3-dimensional graph manifold, a positive answer to the mentioned above question is also given in \cite{NTY21}, and the first result of this paper generalizes the above result in the setting of 3-dimensional graph manifold to a more general setting  with a completely different proof. We emphasize here that the strategy in \cite{NTY21} does not work here.

\begin{thm}
\label{main:thm}
Let $\mathcal{G} = (\Gamma, \{G_v\}, \{G_e\})$ be a graph of groups. Suppose that each edge group is infinite, and each finite height subgroup of vertex group is either 
finite or has finite index.  Let $T$ be the Bass–Serre tree associated to
the decomposition.  Suppose that $G$ acts acylindrically hyperbolic on $T$. Then the following are equivalent.
\begin{enumerate}
    \item $H$ has finite height in $G$.
    \item $(G, T, H)$ is a $A/QI$--triple.
    \item $H$ is strongly quasiconvex and virtually free in $G$.
\end{enumerate}
\end{thm}

\begin{rem}
As mentioned early about the example of torus bundle over circle, the condition that $G$ acts acylindrically on the Bass-Serre tree $T$ in Theorem~\ref{main:thm} is necessary to achieve the above equivalence statements.
\end{rem}

Theorem~\ref{main:thm} could be applied to several classes of graphs of groups, for instance, fundamental groups of 3-dimensional graph manifolds, high dimensional graph manifolds \cite{FLS15} (note that $3$--dimensional graph manifolds are not high dimensional graph manifolds as in the sense of \cite{FLS15}), and  admissible groups \cite{CK02}. We refer the reader to Section~\ref{sec:appl} for further discussion.


In general, a strongly quasiconvex subgroup $P$ of a finitely generated group $G_1$ is not necessary to be strongly quasiconvex in HNN extension $G = G_{1}*_{H_{1}^t = H_2}$ and amalgamation $G = G_1*_{H_1=H_2}G_2$ since $G_1$ could be badly distorted in $G$. It is a natural question to ask in which conditions of $H_1$ and $H_2$, the subgroup $P$ is still strongly quasiconvex in $G$.

Given a finite collection of subgroups $H_1, H_2, \dots, H_n$ in a finitely generated group $G$, we say a finitely generated subgroup $P$ of $G$ has {\it bounded penetration } with respect to the collection $\mathcal{B}$ of left cosets of subgroups $H_i$'s in $G$ if there exists a function $f \colon [0, \infty) \to [0, \infty)$ such that $\diam \bigl ( \mathcal{N}_{r}(P) \cap B) \le f(r)$ for each $B \in \mathcal{B}$.

\begin{thm}
\label{thm:intro2}
Let $G$ be either a HNN extension  $G_{1}*_{H_{1}^t = H_2}$ of a finitely generated group $G_1$ where the stable letter $t$ conjugates the subgroup $H_1$  to the subgroup $H_2$ or $G$ is an amalgamation $G = G_1*_{H_1=H_2}G_2$ such that $H_i$ is undistorted in $G$.
Let $\mathcal{B}$ be the collection of all left cosets of $H_1$ and $H_2$ in $G_1$ in case $G$ is the HNN extension, and let $\mathcal{B}$ be the collection of all left cosets of $H_1$ in $G_1$ in case $G$ is the amalgamation.

Assume that $P \le G_1$ has bounded penetration with respect to the collection $\mathcal{B}$.
Then $P$ is a strongly quasiconvex subgroup in $G_1$ if and only if  $P$ is also a strongly quasiconvex subgroup of $G$.
\end{thm}

\subsection*{Acknowledgments}
We would like to thank Wenyuan Yang for helpful conversations.
 The first author is supported by Project ICRTM04\_2021.07 of the International Centre for Research and Postgraduate Training in
Mathematics, VietNam. The second author was supported by an AMS-Simons Travel Grant.

\section{Preliminary}
\begin{defn}

A subset $A$ of a geodesic space $X$ is called \emph{Morse} if for any $L \ge 1$, $C \ge 0$ there exists $R = R(L,C)$ such that every $(L,C)$--quasi-geodesic in $X$ with endpoints in $A$ is contained in the $R$--neighborhood of $A$.

Let $G$ be a finitely generated group and $H$ a subgroup of $G$.  We say $H$ is \emph{strongly quasiconvex} in $G$ if $H$ is a Morse subset in the Cayley graph $\Gamma(G,S)$ for some (any) finite generating set $S$.
\end{defn}


\begin{defn}
Suppose $(X, d)$ is a Gromov hyperbolic geodesic metric space, and $Y \subset X$. An isometric action $G \curvearrowright X$ is {\it { acylindrical along }}  $Y$ if for every $\epsilon \ge 0$, there are constants $R = R(\epsilon, Y) > 0$ and $M = M(\epsilon, Y) > 0$ so that, whenever $x, y \in Y$ satisfy $d(x, y) \ge R$, \[
\# \bigl \{ g \in G : d(gx, x) \le \epsilon \,\, \text{and} \,\, d(gy, y) \le \epsilon \bigr \} \le M
\]
\end{defn}

The following notion is recently introduced by Abbott-Manning \cite{AM21}.
\begin{defn}
An $A/QI$--triple $(G, X, H)$ consists of a Gromov hyperbolic space $X$ and an isometric action $G \curvearrowright X$ which is acylindrical along the finitely generated infinite subgroup $H$ (i.e, it is acylindrical along along the subspace $Hx$ for some $x \in X$), and so that the orbital map $h \mapsto hx$ gives a quasi-isometric embedding of $H$ into $X$.
\end{defn}

\section{Proof of Theorem~\ref{main:thm} and its application to certain graph of groups}

In this section, we are going to prove Theorem~\ref{main:thm}.
Let $\mathcal{G} = (\Gamma, \{G_v\}, \{G_e\})$ be a graph of groups. We consider a {\it quasi-isometric model} $X$ of $G$ as the following. Choose a maximal tree $\Lambda$ in $\Gamma$. Choose a generating set for $G$ consisting of a finite generating set for each vertex group $G_v$ together with a generator $t_e$ for each edge $e$ of $\Gamma$ not contained in the maximal tree $S$. For edges $e$ in the maximal tree, set $t_e =1$. Then for any edge $e$ and $h \in G_e$, let $h_0, h_1$ be the images of $h$ in the initial and terminal vertex groups of $e$ respectively, then $h_{0} t_e = t_{e} h_1$.
 
 Consider the Bass-Serre tree $T$ whose vertices correspond to cosets of the vertex groups.
For each vertex $v$ in the Bass-Serre tree, take a copy $X_v$ of the Cayley graph of the
corresponding vertex group, but with vertices labelled by elements of the coset. For an edge labelled $e$ in the Bass-Serre tree connecting $v$ to $w$, attach edges between vertices of
the form $gh_{0} \in X_v$ and $g h_{0} t_{e} = g t_{e} h_{1} \in X_w$ for each $h \in G_{e}$. This defines the graph $X$. There is a natural map $\varphi \colon G \to X$ which is an $(L, L)$--quasi-isometric embedding which comes from Normal Form Theorem.


\begin{proof}[Proof of Theorem~\ref{main:thm}]
By Proposition~4.7  and Theorem~1.5 in \cite{AM21}, if $(G, T, H)$ is $A/QI$--triple then $H$ has finite height, and $H$ is stable in $G$, thus the implications $(2) \implies (1)$ and $(3) \implies (1)$ hold. If $H$ is stable in $G$ then $H$ has finite height in $G$ by \cite{Tra19}, thus the implication $(1) \implies (3)$ holds.

For the rest of the proof, we will prove the implication $(1) \implies (2)$. Since $G$ acts acylindrically on the Bass-Serre tree $T$, it follows that $G  \curvearrowright T$ is acylindrical along $H$. 
Hence, to see  $(G, T, H)$ is a $A/QI$--triple, we only need to show that the orbital map $\tau \colon H \to T$ is a quasi-isometric embedding. We will assume that $H$ has infinite index in $G$, otherwise it is obvious that that the implication $(1) \implies (2)$ holds.

First at all, we observe that the intersection of $H$ with any conjugate of a vertex group must be finite. Indeed, by way of contradiction, we assume that $H \cap y G_v y^{-1}$ is infinite for some vertex group $G_v$ and for some $y \in G$. Since $H$ has finite height in $G$, it follows that $H \cap y G_v y^{-1}$ has finite height in $y G_v y^{-1}$. Combining with the fact that $H \cap y G_v y^{-1}$ is infinite, it implies that $H \cap y G_v y^{-1}$ has finite index in $y G_v y^{-1}$ by our assumption. For any $x \in G$, we choose a finite sequence of conjugates of vertex groups $A_0 = yG_v y^{-1}, A_1, \dots, A_m = xA_{0} x^{-1}$ with $A_{i} \cap A_{i+1}$ is infinite with $i \in \{0, 1, \dots, n-1 \}$. Since $H \cap y G_v y^{-1}$ is a finite index subgroup of $ y G_v y^{-1}$ and $y G_v y^{-1} \cap A_1$ is infinite, it follows that $H \cap A_1$ is infinite. Repeating the argument above consecutively for $A_1, A_2, \dots, A_m$ (by replacing $A_0$ by $A_1$), we have that $H \cap A_1$ is a finite index subgroup of $A_1$, $H \cap A_2$ is a finite index subgroup of $A_2$, ...,  $H \cap A_m = H \cap x (yG_vy^{-1}) x^{-1} $ is a finite index subgroup of $x (yG_vy^{-1}) x^{-1}$. Let $s$ be the height of $H$ in $G$. We choose $s+1$ distinct left cosets $g_{0}H, \dots, g_{s} H$ such that the intersection
$g_0 H g_0^{-1} \cap \dots \cap g_s H g_s^{-1}$
is finite. For each $i \in \{1,2, \dots, s \}$, choose $x_i : = g_{i}^{-1}$.  Since $H \cap x_i (yG_vy^{-1}) x_{i}^{-1} $ is a finite index subgroup of $x_{i} (yG_vy^{-1}) x_{i}^{-1}$. It follows that $g_{i} H g_{i}^{-1} \cap yG y^{-1} $ is a finite index subgroup of $y G_v y^{-1}$. We thus have $\bigcap_{i=0}^{s} g_{i} H g_{i}^{-1} \cap y G_v y^{-1}$ is a finite index subgroup of $y G_v y^{-1}$. This is impossible since the intersection
$g_0 H g_0^{-1} \cap \dots \cap g_s H g_s^{-1}$
is finite and  $y G_v y^{-1}$ is infinite.

Let $\mathcal{B}$ be the collection of left cosets of vertex groups $G_{v_1}, G_{v_2}, \dots, G_{v_n}$ in $G$. We define a function $f \colon [0, \infty) \to [0, \infty)$ as follows. For any $r \ge 0$, we define $f(r)$ to be \[
f(r):=\max \bigl \{\diam\bigl(N_r(H)\cap B\bigr): B \in \mathcal{B} \text{ and } B \cap B_{S}(e, r) \neq \emptyset \bigr) \bigr \}
\] 
Note that $f$ is well-defined. Indeed, if $B$ is an element in the collection $\mathcal{B}$ then $B$ is a left coset $gG_{v}$ for some vertex group $G_v$ and for some group element $g$ in $G$.
By Proposition~9.4 in \cite{Hru10} there exists a constant $r' = r'(r, H, gG_v)$ such that 
$
\mathcal N_{r}(H) \cap \mathcal N_{r} (gG_v) \subset \mathcal N_{r'}(H \cap gG_v g^{-1})
$. According to the previous paragraph, the intersection $H \cap gG_v g^{-1}$ is finite, and thus $\mathcal N_{r'}(H \cap gG_v g^{-1})$ is finite. As $\mathcal{N}_{r}(H) \cap B$ is a subset of $\mathcal{N}_{r}(H) \cap \mathcal{N}_{r}(g G_v)$, it follows that $\diam(\mathcal{N}_{r}(H) \cap B)$ is finite.
Since  there are only finitely many left cosets $xA$ that has nonempty intersection with the closed ball $B_{S}(e,r)$, it follows that $f(r)$ is a constant in $[0, \infty)$. 

{\it Claim~1:}
\[
\diam (N_{r}(H) \cap gG_{v} ) \le f(r)
\] for any left coset of a vertex group $G_v$ in $G$.

Indeed, if $d_S(H,B)>r$, then $N_r(H)\cap B$ is empty and then its diameter is zero and less than or equal $f(r)$. We now assume that $d_S(H,B)\leq r$. Then there are $h\in H$, $g\in G$, and $A\in \mathcal{A}$ such that $B=h(gA)$ and $|g|_S\leq r$. By the choice of $f(r)$ the diameter of the set $N_r(H)\cap gA$ is at most $f(r)$. Also $N_r(H)\cap B=h\bigl(N_r(H)\cap gA\bigr)$. Therefore, $\diam \bigl (N_r(H)\cap B \bigr) \le f(r)$. The claim is proved.

Since $H$ has finite height in $G$, it follows from the previous paragraphs that the intersection of $H$ with any conjugate of a vertex group $G_v$ is finite. Since $H \cap g G_{v} g^{-1}$ is finite, it follows that $H$ acts properly on the Bass-Serre tree $T$ and the stabilizer in $H$ of each vertex in $T$ is finite. Hence, it follows from Theorem~7.51 in \cite{DK18} that there exists a finite index subgroup $K$ of $H$ such that $K$ is a free group.

Let $X$ be a quasi-geodesic model space of the graph of groups described in the first paragraph of this section. Let $\varphi \colon G \to X$ be the given $(L, L)$--quasi-isometric embedding.
Let $\{ \gamma_1, \gamma_2, \cdots, \gamma_s \}$ be a finite generating set of $K$. We note that $\gamma_i$ is not conjugate to any vertex group of $G$. Fix a vertex $v_0 \in T$ such that the identity $1 \in X_{v_0}$. For each $i \in \{1, 2, \cdots, s \}$, let $\gamma_{s+i} = \gamma^{-1}_{i}$. Let $\ell_i$ be a geodesic in $X$ connecting $1 = \varphi(1)$ to $\varphi(\gamma_i)$, and let $\ell_{s+i} : = \ell_{i}$. We define
\[
\widetilde{H} : = \bigcup_{i=1}^{2s} \{ g \ell_{i} \, |\, g \in K \}
\]
Then $\widetilde{H}$ is a connected subspace of $X$ and $\varphi(H) \subset \widetilde{H}$. 

{\it {Claim~2:}} There exists a constant $\delta >0$ such that $\diam \bigl ( \widetilde{H} \cap X_v \bigr ) \le \delta$ for each vertex $v \in T$.
Indeed, let $\lambda$ be the maximal length of $\gamma_i$ where $i$ varies from $1$ to $2s$. Let $\delta : = 2 \lambda + L + L f(L \lambda + L)$.
Let $x$ and $y$ be two distinct elements in $\widetilde{H} \cap X_v$. There are elements $g'$ and $g''$  in $K$, and $i, j \in \{1, 2, \cdots, 2s \}$ such that $x \in g' \ell_i$ and $y \in g'' \ell_j$. We thus have $d(x, \varphi(g'\gamma_i)) \le \lambda$ and $d(y, \varphi(g''\gamma_j)) \le \lambda$. Hence $\varphi(g' \gamma_i) \in \mathcal{N}_{\lambda}(X_v) = \mathcal{N}_{\lambda}(\varphi(gG_v)$. Since $\varphi \colon G \to X$ is an $(L, L)$--quasi-isometric, it follows that $g' \gamma_i$ is in  $\mathcal{N}_{L \lambda + L}(gG_v)$. Thus $g' \gamma_i$ belongs to the intersection $ \mathcal{N}_{L \lambda + L}(gG_v) \cap K$. Similarly, we have that the element $g'' \gamma_j$ belongs to the intersection $\mathcal{N}_{L \lambda + L}(gG_v) \cap K$.

By Claim~1, we have $\diam( \mathcal{N}_{L\lambda + L}(K) \cap \mathcal{N}_{L \lambda + L}(gG_v) \le f(L \lambda + L)$. Thus, $d_{G}(g' \gamma_i,g'' \gamma_{j}) \le f(L \lambda + L) $. It follows that
\begin{align*}
    d(x,y) &\le 2\lambda + d(\varphi(g' \gamma_i), \varphi(g'' \gamma_j) \\
    &\le 2 \lambda + L d_{G}(g' \gamma_i,g'' \gamma_{j}) + l \\
    &\le 2 \lambda + L + L\, f(L \lambda + L) = \delta
\end{align*}
Therefore, $\diam( \widetilde{H} \cap X_v) \le \delta$, Claim~2 is confirmed.

Let $e_1 \cdot e_2 \cdots e_s$ be the geodesic edge path in the Bass-Serre tree $T$ connecting $v_0$ to $h(v_0)$. Let $v_i$ be the endpoint of $e_i$. Let $\beta$ be a path in $\widetilde H$ connecting $\varphi(1) =1$ to $\varphi(h)$. We note that the path $\beta$ must pass through vertex and edge spaces $X_{v_{0}}$, $X_{e_{1}}$, $X_{v_{1}} \cdots, X_{v_{s}}$.
Let $x_0 \in \beta \cap X_{v_0} \cap X_{e_1} $ be the last point of $\beta$ which entering $X_{e_1}$. Note that there is an edge in $X_{e_1}$ connecting $x_0$ to a point in $X_{v_1}$, we denote this point by $y_0$ (so $d(x_0, y_0) =1$). Similarly, we define points $x_1, y_1, \cdots, x_{s-1}, y_{s-1}$ in a similar manner. We let $x_{s} : = \varphi(h)$. By Claim~2 we have $d(1, x_0) \le \delta$, $d(y_i, x_{i+1}) \le \delta$ for each $0 \le i \le s-1$. We also have $d(x_i, y_i) = 1$ for $0 \le i \le s-1$.

Let $\Lambda > 0$ be a constant  such that $\widetilde H  \subset \mathcal{N}_{\Lambda} (\varphi(K))$. It follows that for each $0 \le i \le s$, there exists $h_i \in H$ such that $d(x_i, \varphi(h_i))$ is bounded above by $ \Lambda$ (when $i =s$, we choose $h_s = h$).
Let $\hat{\epsilon} := \max \{ |h|_{H} \, : \, h \in B_{G}(1, L + L(1 + \delta + 2\Lambda) \}$.

By triangle inequality, we have $ d(\varphi(1), \varphi(h_0)) \le d(\varphi(1), x_0) + d(x_0, \varphi(h_0)) 
    \le \delta + \Lambda$. Since the map $\varphi$ is a $(L, L)$--quasi-isometric embedding map, we have that  $d_{G}(1, h_0) \le L d(\varphi(1), \varphi(h_0)) + L \le L(\delta + \Lambda) + L$. It follows that $d_{H}(1, h_0)$ is bounded above by $ \hat{\epsilon}$.
    For each $0 \le i \le s$, we have 
$
    d(\varphi(h_i), \varphi(h_{i+1}) \le d(\varphi(h_i), x_i) + d(x_i,y_{i}) + d(y_i, x_{i+1}) + d(x_{i+1}, \varphi(h_{i+1}) 
    \le 2\Lambda + \delta +1
$
Hence $d_{G}(h_i, h_{i+1}) \le L( 2 \Lambda + \delta+1) + L$. It follows that, $d_{H}(h_{i}, h_{i+1}) \le \hat{\epsilon}$.
As a consequence, 
$d_{H}(1, h) \le d_{H}(1, h_0) + \sum_{i=0}^{s-1} d_{H}(h_i, h_{i+1}) \le (s+1) \hat{\epsilon} \le \hat{\epsilon} + \hat{\epsilon} d_{T}(v_0, h(v_0))$.

On the other hand,  the orbital map of any isometric action is Lipschitz (see Lemma~I.8.18 in \cite{BH99}) , and thus  the orbital map $\tau \colon H \to T$ is a Lipschitz map.

Therefore, we can conclude that the orbital map $\tau$ is a quasi-isometric embedding. The implication $(1) \implies (2)$ is verified.
\end{proof}

\subsection{Application of Theorem~\ref{main:thm} to certain graph of groups}
\label{sec:appl}
 In this section, we are going give some examples of graph of groups which satisfy hypotheses of Theorem~\ref{main:thm}. Therefore, the conclusion of Theorem~\ref{main:thm} holds for these classes of graph of groups.
 
 \begin{lem}
\label{central}
Let $G$ be a group such that the centralizer $Z(G)$ of $G$ is infinite. Let $H$ be a finite height infinite subgroup of $G$. Then $H$ must have finite index in $G$.
\end{lem}
\begin{proof}
We first assume that $Z(G) \cap H$ has infinite index in $Z(G)$. Then there is an infinite sequence $(t_n)$ of elements in $Z(G)$ such that $t_i(Z(G) \cap H)\neq t_j(Z(G) \cap H)$ for $i\neq j$. Therefore, it is straightforward that $t_iH\neq t_jH$ for $i\neq j$. Also, $\bigcap t_nHt_n^{-1}=H$ is infinite. This contradicts to the fact that $H$ has finite height. Therefore, $Z(G) \cap H$ has finite index in $Z(G)$, In particular, $Z(G) \cap H$ is infinite. Assume that $H$ has infinite index in $G$. Then there is an infinite sequence $\{g_nH\}$ of distinct left cosets of $H$. However, $\bigcap g_nHg_n^{-1}$ is infinite since it contains the infinite subgroup $Z(G) \cap H$. This contradicts to the fact that $H$ has finite height. Therefore, $H$ must have finite index in $G$.
\end{proof}

\subsection*{3-dimensional graph manifolds}
We revisit a result of \cite{NTY21} in the setting of 3-dimensional graph manifold.
Let $M$ be a compact, connected, orientable, irreducible 3-manifold with empty or toroidal boundary. $M$ is called \emph{geometric} if its interior admits a geometric structure in the sense of Thurston. Such structures are $S^3$, $\mathbb{E}^3$, $\mathbb{H}^3$, $S^{2} \times \mathbb{R}$, $\mathbb{H}^{2} \times \mathbb{R}$, $\widetilde{SL(2, \mathbb{R})}$, Nil, and Sol. If $M$ is not geometric, then $M$ is called a {\it nongeometric $3$--manifold}. By geometric decomposition  of $3$--manifolds, there is a nonempty minimal union $\mathcal{T} \subset M$ of disjoint essential tori and Klein bottles, unique up to isotopy, such that each component of $M \backslash \mathcal{T}$ is either a Seifert fibered piece or a hyperbolic piece.
The manifold $M$ is called a {\it graph manifold} if all the pieces $M_1$, $M_2, \dots, M_k$ of $M \backslash \mathcal{T}$ are Seifert manifolds. Passing to a finite cover, we can assume that the base orbifold of each $M_i$ is orientable and hyperbolic The fundamental group $\pi_1(M)$ has the graph of groups structure where it has one vertex labelled with $\pi_1(M_i)$ for each $i$, and each edge labelled by $\mathbb{Z}^2$ for each JSJ torus. Let $T$ be the Bass-Serre tree associated to the graph of groups structure of $\pi_1(M)$. Then the action of $\pi_1(M)$ on $T$ is acylindrical (see Theorem~7.27 in \cite{ABO19}). Each vertex group $\pi_1(M_i)$ of $\pi_1(M)$ has centralizer $\mathbb{Z}$ which is infinite, and thus by Lemma~\ref{central}, each finite height subgroup in each vertex group of $\pi_1(M)$ is either finite or has finite index in $\pi_1(M_i)$. As a consequence, Theorem~\ref{main:thm} can be applied into graph 3-manifolds. 

\subsection*{High dimensional graph manifolds}
In \cite{FLS15}, Frigerio-Lafont- Sisto study a particular class of graph of groups which they called high dimensional graph manifold which are obtained by gluing together pieces of the form $T_{i}^{k} \times N_i$ where $T_{i}^{k}$ is a $k$--dimensional torus, and $N_i$ is the manifold compactification of a complete, finite volume hyperbolic manifold $N_i$ of dimension $\ge 3$.
Such manifolds are called \emph{irreducible} if the fibrations on adjacent pieces are ``transverse''. We refer the reader to \cite{FLS15} for precise definitions. We emphasize here that $3$--dimensional graph manifolds are not high dimensional graph manifolds as in the sense of \cite{FLS15}.

When $M$ is an irreducible high dimensional graph manifold then $\pi_1(M)$ is a graph of groups where each vertex groups, edge groups are undistorted in $\pi_1(M)$ (see Corollary~7.13 in \cite{FLS15}). Also $\pi_1(M)$ acts acylindrically on its Bass-Serre tree $T$ (see Proposition~6.4 in \cite{FLS15}). The fundamental group of each piece $M_i$ of a high dimensional graph manifold is the product of $\pi_1(N_i)$ with $\mathbb{Z}^k$. Hence the centralizer $Z(\pi_1(M_i))$ is infinite. It follows from Lemma~\ref{central} that each finite height subgroup in each vertex group of $\pi_1(M)$ is either finite or has finite index in $\pi_1(M_i)$. Therefore the fundamental groups of high dimensional graph manifolds satisfy the hypotheses of Theorem~\ref{main:thm}. As a result, finite height, stability and  $A/QI$--triple of subgroups in high dimensional graph manifolds are equivalent. 

\subsection*{Admissible groups} 
This class of groups firstly introduced by Croke-Kleiner in \cite{CK02}. This is a particular class of graph of groups that includes fundamental groups of $3$--dimensional graph manifolds and torus complexes. The admissible group  is modeling on the JSJ structure of graph manifolds where the Seifert fibration is replaced by the following central extension of a general hyperbolic group:
\begin{equation}\label{centralExtEQ}
1\to Z(G_v)=\mathbb Z\to G_v\to H_v\to 1    
\end{equation}
Precisely, a graph of groups $\mathcal{G}$ is \emph{admissible} if  $\mathcal{G}$ is a finite graph with at least one edge, each vertex group ${ G}_v$ has center $Z({ G}_v) \cong \Z$, ${ H}_v \colon = { G}_{v} / Z({ G}_v)$ is a non-elementary hyperbolic group, and every edge subgroup ${ G}_{e}$ is isomorphic to $\Z^2$. Moreover, 
\begin{enumerate}
    \item Let $e_1$ and $e_2$ be distinct directed edges entering a vertex $v$, and for $i = 1,2$, let $K_i \subset { G}_v$ be the image of the edge homomorphism ${G}_{e_i} \to {G}_v$. Then for every $g \in { G}_v$, $gK_{1}g^{-1}$ is not commensurable with $K_2$, and for every $g \in  G_v \setminus K_i$, $gK_ig^{-1}$ is not commensurable with $K_i$.
    \item For every edge group ${ G}_e$, if $\alpha_i \colon { G}_{e} \to { G}_{v_i}$ is the edge monomorphism, then the subgroup generated by $\alpha_{1}^{-1}(Z({ G}_{v_1}))$ and $\alpha_{2}^{-1}(Z({ G}_{v_1}))$ has finite index in ${ G}_e$.
\end{enumerate}
A group $G$ is \emph{admissible} if it is the fundamental group of an admissible graph of groups.

Since the centralizer of each vertex group is infinite, it follows from Lemma~\ref{central} that each finite height subgroup in each vertex group  is either finite or has finite index in the vertex group. In addition, $G$ acts acylindrically on its Bass-Serre tree (see part~(1) of Lemma~3.4 in \cite{CK02}). Therefore Theorem~1.2 applies to admissible groups. We note that when we restrict admissible groups to be CAT(0), such a result has been addressed in \cite{NY20}. 

\section{Morse subsets in group amalgamations and HNN extensions}
In this section, we are going to prove Theorem~\ref{thm:intro2}.
\begin{defn}
Let $X$ be a geodesic space and let $A$ be a subset of $X$. Let $\mathcal{B}$ be a collection of subspaces of $X$. Let $f\colon [0,\infty)\rightarrow [0,\infty)$ be a function. We say $A$ has {\it $f$--bounded penetration} with respect to $\mathcal{B}$ if for each $B$ in $\mathcal{B}$ and each positive number $r$ we have $\diam( N_r(A)\cap B ) \le f(r)$. We say $A$ has \emph{bounded penetration} with respect to $\mathcal{B}$ if  $A$ has \emph{$f$--bounded penetration} with respect to $\mathcal{B}$ for some function $f$.
\end{defn}

\begin{lem}
\label{le1}
Let $\lambda \geq 1$, $\epsilon \geq 0$, and $K>2$ be constants. There are two numbers $\lambda_1=\lambda_1(\lambda,\epsilon, K)\geq 1$ and $\epsilon_1=\epsilon_1(\lambda,\epsilon, K)\geq 0$ such that the following holds. Let $X$ be a geodesic space and let $\gamma$ be a continuous $(\lambda,\epsilon)$--quasigeodesic parametrized by arc its length in $X$. Let $x_1$ and $x_2$ be two points in $X$ and let $x'_i$ be a nearest point projection of each $x_i$ onto $\gamma$. For each $i$ let $\alpha_i$ be a geodesic in $X$ connecting $x_i$ and $x'_i$ and let $\gamma_1$ be the subpath of $\gamma$ connecting $x'_1$ and $x'_2$. Assume that $d(x'_1,x'_2)\geq K \max\{d(x_1,x'_1),d(x_2,x'_2)\}$. Then the concatenation $\beta=\alpha_1\gamma_1\alpha_2$ is a continuous $(\lambda_1,\epsilon_1)$--quasigeodesic parametrized by its arc length.
\end{lem}

\begin{proof}
For each two points $x$ and $y$ in $\beta$ we denote $\ell_{\beta}(x,y)$ be the length of the subpath of $\beta$ connecting $x$ and $y$. We will prove that for each two points $u$ and $v$ in $\beta$ we have $\ell_{\beta}(u,v)\leq \lambda_1 d(u,v)+\epsilon_1$ where $\lambda_1$ and $\epsilon_1$ only depend on $\lambda$, $\epsilon$, and $K$. Since the paths $\alpha_1$, $\alpha_2$, and $\gamma$ are all $(\lambda,\epsilon)$--quasigeodesic, the inequality is obvious if $u$ and $v$ both lie in one of these path. Therefore, we can assume that $u$ and $v$ lie in two different paths among $\alpha_1$, $\alpha_2$, and $\gamma_1$. Without the loss of generality we can assume only two cases occur:

\textbf{Case 1:} $u$ lies in $\alpha_1$ and $v$ lies in $\gamma_1$. Since $x'_1$ is a nearest point projection of $u$ of $\gamma$, we have $d(u,v)\geq d(u,x'_1)\geq d(x'_1,v)-d(u,v).$ Therefore,
$
d(u,v)\geq \frac{1}{2}d(x'_1,v)\geq \frac{1}{2}\biggl(\frac{\ell_\beta(x'_1,v)}{\lambda}-\frac{\epsilon}{\lambda}\biggr)$.
This implies that $\ell_\beta(x'_1,v)\leq 2\lambda d(u,v)+\epsilon$. Also, $\ell_\beta(u,x'_1)=d(u,x'_1)\leq d(u,v)$, thus, $\ell_\beta(u,v)=\ell_\beta(u,x'_1)+\ell_\beta(x'_1,v)\leq (2\lambda+1) d(u,v)+\epsilon$.

\textbf{Case 2:} $u$ lies in $\alpha_1$ and $v$ lies in $\alpha_2$. Then
\begin{align*}
    d(u,v)&\geq d(x'_1,x'_2)-d(x'_1,u)-d(x'_2,v)\\&\geq d(x'_1,x'_2)-d(x'_1,x_1)-d(x'_2,x_2) \\&\geq d(x'_1,x'_2)-\frac{1}{K}d(x'_1,x'_2)-\frac{1}{K}d(x'_1,x'_2)\\&\geq \biggl(\frac{K-2}{K}\biggr) d(x'_1,x'_2) \geq \biggl(\frac{K-2}{K}\biggr)\biggl(\frac{\ell_\beta(x'_1,x'_2)}{\lambda}-\frac{\epsilon}{\lambda}\biggr).
\end{align*}
This implies that $\ell_\beta(x'_1,x'_2)\leq \biggl(\frac{K\lambda}{K-2}\biggr)d(u,v)+\epsilon$.
Also, $\ell_\beta(u,x'_1)=d(u,x'_1)\leq d(x_1,x'_1)\leq d(x'_1,x'_2)\leq \ell_\beta(x'_1,x'_2)$.
Similarly, $\ell_\beta(x'_2,v)\leq \ell_\beta(x'_1,x'_2)$.
Therefore, $\ell_\beta(u,v) = \ell_\beta(u,x'_1) + \ell_\beta(x'_1,x'_2)+\ell_\beta(x'_2,v) \leq 3\ell_\beta(x'_1,x'_2)\leq  \biggl(\frac{3K\lambda}{K-2}\biggr)d(u,v)+3\epsilon$
\end{proof}

\begin{lem}[Lemma 3.3 in \cite{Tra19}]
\label{ll2}
For each $C>1$ and $\rho \in (0,1]$ there is a constant $L=L(C,\rho)\geq 1$ such that the following holds. Let $r$ be an arbitrary positive number and $\gamma$ a continuous path with the length less than $Cr$. Assume the distance between two endpoints $x$ and $y$ is at least $r$. Then there is an $(L,0)$--quasi-geodesic $\alpha$ connecting two points $x$, $y$ such that the image of $\alpha$ lies in the $\rho r$--neighborhood of $\gamma$ and the length of $\alpha$ is less than or equal to the length of $\gamma$.
\end{lem}

\begin{prop}
\label{cool}
Let $\lambda \geq 1, \epsilon \geq 0$, and $\sigma \geq 0$ be constants. Let $M\colon [1,\infty)\times [0,\infty)\rightarrow [0,\infty)$ be a Morse gauge and $f\colon [0,\infty)\rightarrow [0,\infty)$ be a function. Then there is a number $b=b(\lambda,\epsilon,\sigma,M)$ such that the following holds. Let $X$ be a geodesic space and let $A$ be an $M$--Morse subset of $X$. Assume that $A$ has $f$--bounded penetration with respect to a collection $\mathcal{B}$ of $(\lambda,\epsilon,\sigma)$--quasiconvex subspaces of $X$. Then for each $B$ in $\mathcal{B}$ and each positive number $r$, we have $\diam \bigl ( N_r(A) \cap B \bigr ) \le 7r +b$.
\end{prop}

\begin{proof}
By Lemma 1.11 of \cite{BH99} III.H there are $\lambda_0,\epsilon_0$, and $\epsilon_0$ depending only on $\lambda,\epsilon$, and $\epsilon$ such that any two points of a subset $B$ in $\mathcal{B}$ are joined by a continuous $(\lambda_0,\epsilon_0)$--quasigeodesic parametrized by its length that lies in the $\sigma_0$--neighborhood of $B$. Let $\lambda_1=\lambda_1(\lambda_0,\epsilon_0,3)\geq 1$ and $\epsilon_1=\epsilon_1(\lambda_0,\epsilon_0,3)\geq 0$ be constants in Lemma~\ref{le1}. Let $C=M(\lambda_1,\epsilon_1)+\sigma_0$ and let $b=f\bigl(f(C)+\sigma_0\bigr)$. 

Let $B$ be a subspace in $\mathcal{B}$ and let $r$ be a positive number. We will prove that $\diam \bigl ( N_r(A) \cap B \bigr ) \le 7r +b$. If $r\leq f(C)+\sigma_0$, then $\diam \bigl ( N_r(A)\cap B \bigr ) \le b:=f\bigl(f(C)+\sigma_0\bigr) \le 7r +b $. We now assume that $r>f(C)+\sigma_0$.

Suppose by way of contradiction that $\diam \bigl ( N_r(A)\cap B \bigr ) > 7r +b$. Let $z_1$ and $z_2$ be two points in $N_r(A)\cap B$ such that  $d(Z_1, Z_2)  > 7r+b$. For each $i$, let $x_i$ be a point in $A$ such that $d(z_i,x_i)<r$. We connect $z_1$ and $z_2$ by a continuous $(\lambda_0,\epsilon_0)$--quasigeodesic $\gamma$ parametrized by its length that lies in the $\sigma_0$--neighborhood of $B$. For each $i$ let $x'_i$ be a nearest point projection of $x_i$ on $\gamma$. Then for each $i$ we have $d(x_i,x'_i)\leq r$ and which implies that $d(z_i,x'_i)\leq 2r$. Therefore,
$$d(x'_1,x'_2)\geq d(z_1,z_2)-d(z_1,x'_1)-d(z_2,x'_2)\geq 7r-2r-2r\geq 3r.$$
This implies that $d(x'_1,x'_2)\geq 3\max\{d(x_1,x'_1), d(x_2,x'_2)\}$.

For each $i$ let $\alpha_i$ be a geodesic connecting $x_i$ and $x'_i$ and let $\gamma_1$ be the subpath of $\gamma$ connecting $x'_1$ and $x'_2$. Then the concatenation $\beta=\alpha_1\gamma_1\alpha_2$ is a $(\lambda_1,\epsilon_1)$--quasigeodesic by the choice of $\lambda_1$ and $\sigma_1$. Therefore, $\beta$ lies in the $M(\lambda_1,\epsilon_1)$--neighborhood of $A$. In particular, each $x'_i$ lies in the $M(\lambda_1,\epsilon_1)$--neighborhood of $A$. Since each $x'_i$ also lies in the $\sigma_0$--neighborhood of $B$, we can chose $b_i$ in $B$ such that $d(x'_i,b_i)<\sigma_0$. Therefore, each $b_i$ lies in the the subset $N_C(A)\cap B$ by the choice of $C$ and $$d(b_1,b_2)\geq d(x'_1,x'_2)-d(x'_1,b_1)-d(x'_2,b_2)\geq 3r-2\sigma_0> f(C).$$ This implies that $\diam\bigl(N_r(A)\cap B\bigr)> f(C)$ which is a contradiction. Therefore, the diameter of the subset $N_r(A)\cap B$ is at most $7r+b$.  
\end{proof}

\begin{lem}
\label{l1}
Let $G_1 \le G$ be the groups given by Theorem~\ref{thm:intro2}. Then $G_1$ is undistorted in $G$.
\end{lem}

\begin{proof}
Let $S_H$ be a finite generating set of $H$. Then $S_K=\{\phi(h)|h\in H\}$ is a finite generating set of $K$. Let $T$ be a finite generating set of $G_1$ that contains $S_H\cup S_K$. Then $S=T\cup \{t\}$ be a finite generating set of the HNN extension $G=\langle G_1, t| t^{-1}ht=\phi(h), h\in H\rangle$. With those chosen finite generating sets, we can assume that Cayley graphs $\Gamma(H,S_H)$ and $\Gamma(K,S_K)$ are subgraphs of Cayley graph $\Gamma(G_1,S)$ and Cayley graph $\Gamma(G_1,S)$ is a subgraph of Cayley graph $\Gamma(G,T)$.

Since $H$ and $K$ are undistorted in $G$, there is a positive integer $M$ such that
$|h|_{S_H}\leq M |h|_T \text{ for each $h\in H$}$
and 
$|k|_{S_K}\leq M |k|_T \text{ for each $k\in K$}$
We will prove that $|g|_{S}\leq M |g|_T \text{ for each $g\in G_1$}$.

In fact, let $\alpha$ be a geodesic in $\Gamma(G,T)$ connecting the identity $e$ and $g$. Then $\alpha$ is decomposed as $\alpha_0\beta_1\alpha_1\cdots \beta_n\alpha_n$ such that each subpath $\alpha_i$ lies completely in $\Gamma(G_1,S)$ and each subpath $\beta_i$ intersects to $\Gamma(G_1,S)$ only at its endpoints $x_i$ and $y_i$. We observe that each $x_i^{-1}y_i$ is a group element in $H\cup K$. Therefore, we can connect $x_i$ and $y_i$ by a path $\beta'_i$ in $\Gamma(G_1,S)$ such that $\ell(\beta'_i)\leq M\ell(\beta_i)$. Therefore, the path $\alpha'=\alpha_0\beta'_1\alpha_1\cdots \beta'_n\alpha_n$ lies completely in $\Gamma(G_1,S)$ and $\ell(\alpha')\leq M \ell(\alpha)$. Therefore, $|g|_{S}\leq M |g|_T$. This implies that $G_1$ is also undistorted in $G$.
 \end{proof}
 

We now ready to give the proof of Theorem~\ref{thm:intro2}.
\begin{proof}[Proof of Theorem~\ref{thm:intro2}]

Suppose that $P$ is strongly quasiconvex in $G$. Since $P \le G_1 \le G$ and $G_1$ is undistorted in $G$ by Lemma~\ref{l1}, it follows from Proposition~4.10 in \cite{Tra19} that $P$ is strongly quasiconvex in $G_1$.

For the rest of the proof, we are going to verify that $P$ is strongly quasiconvex in $G$ provided $P$ is strongly quasiconvex in $G_1$. We will give the proof for the case HNN extension. The case of amalgamated free products is proved similarly.

Let $S_H$ be a finite generating set of $H$. Then $S_K= \{\phi(h) \,\, \bigl |\,\, h\in H\}$ is a finite generating set of $K$. Let $S$ be a finite generating set of $G_1$ that contains $S_H\cup S_K$. Then $T=S\cup \{t\}$ is a finite generating set of the HNN extension $G= \bigl \langle G_1, t\,|\, t^{-1}ht=\phi(h), h\in H \bigr \rangle$. With those chosen finite generating sets, we can assume that Cayley graphs $\Gamma(H,S_H)$ and $\Gamma(K,S_K)$ are subgraphs of Cayley graph $\Gamma(G_1,S)$ and Cayley graph $\Gamma(G_1,S)$ is a subgraph of Cayley graph $\Gamma(G,T)$.

Let $M$ be a Morse gauge of $P$ in Cayley graph $\Gamma(G_1,S)$. By Lemma~\ref{l1} the subgroup $G_1$ is undistorted in $G$. Thus, there are $\lambda_0 \geq 1$ and $\sigma_0\geq 0$ such that each two points $x$ and $y$ in $G_1$ we have $d_T(x,y)\leq d_S(x,y)\leq \lambda_0 d_T(x,y)+\epsilon_0$.   Since two subgroups $H$ and $K$ are undistorted in $G$, they are also undistorted in $G_1$. By Proposition~\ref{cool} there is a constant $b$ such that for each positive integer $r$ and each left coset $B$ in $\mathcal{B}$ the diameter of the intersection between the $r$--neighborhood of $P$ and $B$ is at most $7r+b$ with respect to the metric $d_S$. . 

To see that $P$ is a Morse subset of $G$, we let $\gamma:[u,v]\to \Gamma(G,T)$ be a $(\lambda,\epsilon)$--quasigeodesic that connects two points in $P$. We are going to show that $\gamma \subset \mathcal{N}_{r}(P)$ for some constant $r$ depending only on   $\lambda, \lambda_0,
\epsilon,\epsilon_0,b$ and $M$.

By Lemma Lemma 1.11 of \cite{CS15} III.H we can assume that $\gamma$ is continuous and parametrized by its arc length. Let $r=\max\{d_T(x,P)|x\in \gamma\}$ and let $t_0 \in [u,v]$ such that $d(\gamma(t_0),P) = r $. Let $L= 72\lambda (\epsilon+b+1)(\lambda_0+\epsilon_0)$. Assume that $r> L$ and we will prove that $r$ is less than some number only depending on $\lambda, \lambda_0,
\epsilon,\epsilon_0,b$ and $M$.

Let $[t_1,t_2]$ be a subinterval containing $t_0$ in $[a,b]$ such that $d_T(\gamma(t_i),P)=r/L$ for each $i$ and the image of the path $\gamma_1=\gamma_{|[t_1,t_2]}$ lies outside the $r/L$--neighborhood of $P$. We observe that $|t_i-t_0|\geq r-r/L\geq r/2$. Therefore, $|t_2-t_1|=|t_2-t_0|+|t_1-t_0|\geq r$. Let $K=6\lambda (\epsilon+b+9)(\lambda_0+\epsilon_0)$ We now consider two cases:

\textbf{Case 1:} $|t_2-t_1|<Kr$

For each $i$ let $\alpha_i$ be a geodesic in $\Gamma(G,T)$ with length $r/L$ connecting $\gamma(t_i)$ to a point $x_i$ in $P$. Let $\alpha$ be the concatenation $\alpha_1 \gamma_1 \alpha_2$. We have that
$\ell(\alpha)\leq \ell(\alpha_1)+\ell(\gamma_1)+\ell(\alpha_2)\leq (K+2)r$ 
Moreover, \begin{align*}
    d_S(x_1,x_2)\geq d_T(x_1,x_2)&\geq d_T(\gamma(t_1),\gamma(t_2))-2r/L\\&\geq |t_2-t_1|/\lambda-\epsilon-2r/L\\
  &\geq r/\lambda-\epsilon-2r/L\\&\geq\frac{r}{2\lambda}-\frac{2r}{L}\geq \frac{r}{4\lambda}.
\end{align*}

We now construct a continuous path $\alpha'$ in $\Gamma(G_1,S)$ connecting two points $x_1$ and $x_2$ such that $\alpha'\cap \gamma_1$ is non-empty and $\ell(\alpha')\leq (\lambda_0+\epsilon_0)(K+2)r$.
If $\alpha$ lies completely inside $\Gamma(G_1,S)$ then $\alpha'=\alpha$ is a desired path. We now assume that $\alpha$ does not lie completely inside $\Gamma(G_1,S)$. Then $\alpha$ is decomposed as $\sigma_0\eta_1\sigma_1\cdots \eta_n\sigma_n$ such that the following holds:
\begin{enumerate}
    \item Each subpath $\sigma_i$ lies completely in $\Gamma(G_1,S)$; and
    \item Each subpath $\eta_i$ is not degenerate and it intersects to $\Gamma(G_1,S)$ only at its endpoints $y_i$ and $z_i$.
\end{enumerate}
We observe that each $y_i^{-1}z_i$ is a group element in $H\cup K$. Therefore, we can connect $y_i$ and $z_i$ by a continuous $(\lambda_0,\epsilon_0)$ quasigeodesic (parametrized by arc length) $\eta_i'$ of $\Gamma(G,T)$ which lies completely inside $\Gamma(G_1,S)$. In particular, $\ell(\eta_i')\leq (\lambda_0+\epsilon_0)d_S(y_i,z_i)\leq (\lambda_0+\epsilon_0)\ell(\eta_i)$.
Let $\alpha'=\sigma_0\eta'_1\sigma_1\cdots \eta'_n\sigma_n$. Then
$\ell(\alpha')\leq (\lambda_0+\epsilon_0)\ell(\alpha)\leq (\lambda_0+\epsilon_0)(K+2)r$.

Suppose by way of contradiction that $\alpha'\cap \gamma_1$ is empty. This means that $\gamma_1$ is a subpath of some $\eta_i$ that is replaced by $\eta'_i$ in $\alpha'$.  Therefore, there is some $i$ such that $y_i$ lies in $\alpha_0$ and $z_i$ lies in $\alpha_1$. This implies that 
\begin{align*}
    d_S(y_i,z_i)\geq d_T(y_i,z_i)&\geq d_T(\gamma(t_1),\gamma(t_2))-2r/L\\&\geq |t_2-t_1|/\lambda-\epsilon-2r/L\\&\geq r/\lambda-\epsilon-2r/L\\&\geq\frac{r}{2\lambda}-\frac{2r}{L}\\&\geq \frac{r}{4\lambda}\geq \frac{r}{8\lambda}+b>\frac{7r(\lambda_0+\epsilon_0)}{L}+b.
\end{align*}

Since $y_i, z_i \in \Gamma(G_1, S)$ and $y^{-1}_{i} z_{i} \in H \cup K$, it follows that $y_i$ and $z_i$ both lies in some left coset $B$ in $\mathcal{B}$. Therefore, $y_i$ and $z_i$ both lies in the intersection $N_{r/L}(P)\cap B$ where $N_{r/L}(P)$ is the $(r/L)$--neighborhood of $P$ in $\Gamma(G,T)$. Also the distance between two distinct points in $G_1$ with respect to the metric $d_S$ is at most $(\lambda_0+\epsilon_0)$ times the distance between them with respect to the metric $d_T$. This implies that  $y_i$ and $z_i$ both lie in the intersection $N^S_{r(\lambda_0+\epsilon_0)/L}(P)\cap B$ where $N^S_{r(\lambda_0+\epsilon_0)/L}(P)$ is the $(r/L)$--neighborhood of $P$ in $\Gamma(G_1,S)$. Therefore by Proposition~\ref{cool} we have $d_S(y_i,z_i)\leq 7r(\lambda_0+\epsilon_0)/L+b$ which is a contradiction. Therefore, $\alpha'\cap \gamma_1$ is non-empty. This implies that $\alpha'$ contains a point outside the $(r/L)$--neighborhood of $P$.

We recall that $\alpha'\subset \Gamma(G_1,S)$ has length less than $(\lambda_0+\epsilon_0)(K+2)r$ and connects two points $x_1$ and $x_2$ in $P$ with $d_S(x_1,x_2)\geq r/(4\lambda)$. Let $r_1=r/(4\lambda)$, $C=4\lambda(\lambda_0+\epsilon_0)(K+2)>1$, $\rho=2\lambda/L \in (0,1]$. Then $\rho r_1=r/(2L)$. Moreover, the length of $\alpha'$ is less than $Cr_1$ and the distance between the two endpoints $x_1$ and $x_2$ of $\alpha'$ is greater than $r/(4\lambda)=r_1$. By Lemma~\ref{ll2} there are a number $\lambda_1=\lambda_1(C,\rho)\geq 1$ and a $(\lambda_1,0)$--quasigeodesic $\alpha''$ in $\Gamma(G_1,S)$ connecting two points $x_1$ and $x_2$ such that $\alpha''$ lies in the $(r/(2L))$--neighborhood of $\alpha'$ with respect to the metric $d_S$. Since the distance between any two points in $G_1$ with respect to the metric $d_S$ is greater than or equal the the distance between these two points with respect to the metric $d_T$. Thus, $\alpha''$ lies in the $(r/(2L))$--neighborhood of $\alpha'$ with respect to the metric $d_T$. Therefore, $\alpha''$ contains a points outside the $(r/(2L))$--neighborhood of $P$ with respect to the metric $d_T$. Also, $\alpha''$ lies in the $M(\lambda_1,0)$--neighborhood of $P$ with respect to the metric $d_S$ (therefore also with respect to the metric $d_T$) since $P$ is an $M$--Morse subset in $\Gamma(G_1,S)$. Therefore, $r\leq 2LM(\lambda_1,0)$.

\textbf{Case 2:} $|t_2-t_1|\geq Kr$. 

Let $t'_1,t'_2$ be subinterval of $[t_1,t_2]$ such that $|t'_2-t'_1|=Kr$. Then the image of the path $\gamma_2=\gamma_{|[t'_1,t'_2]}$ lies outside the $(r/L)$--neighborhood of $P$.
 For each $i$ let $\beta_i$ be a geodesic in $\Gamma(G,T)$ with at most $r$ connecting $\gamma(t'_i)$ to a point $x'_i$ in $P$. Let $\beta$ be the concatenation $\beta_1 \gamma_2 \beta_2$. Therefore, the length of $\beta$ is at most $(K+2)r$. Moreover, \begin{align*}
    d_S(x'_1,x'_2)\geq d_T(x'_1,x'_2)&\geq d_T(\gamma(t'_1),\gamma(t'_2))-2r\\&\geq |t'_2-t'_1|/\lambda-\epsilon-2r\\&\geq (Kr)/\lambda-\epsilon-2r\\&\geq\frac{Kr}{2\lambda}-2r\geq r.
\end{align*}
Using an analogous argument as in Case 1, we can construct a continuous path $\beta'$ in $\Gamma(G_1,S)$ connecting two points $x'_1$ and $x'_2$ such that $\beta'\cap \gamma_2$ is non-empty and $\ell(\beta')\leq (\lambda_0+\epsilon_0)(K+2)r$. Therefore, $\beta'$ contains a point outside the $(r/L)$--neighborhood of $P$. Using an analogous argument as in Case 1 again, there is a number $\lambda_2$ depending only on $\lambda_0, \epsilon_0,K, \lambda$ and $L$ and a $(\lambda_2,0)$--quasigeodesic $\beta''$ in $\Gamma(G_1,S)$ connecting two points $x'_1$ and $x'_2$ such that $\beta''$ contains a point outside the $(r/(2L))$--neighborhood of $P$. Therefore, $r\leq 2LM(\lambda_2,0)$.
\end{proof}

\bibliographystyle{alpha}


\end{document}